\documentclass[12pt]{article}
\usepackage{amssymb}
\usepackage{amsmath}

\begin{document}

\begin{center}
\textbf{\large Thickness conditions and Littlewood--Paley sets}
\footnote{This study was carried out within The National Research
University Higher School of Economics' Academic Fund Program in
2013-2014, research grant No. 12-01-0079.}
\end{center}

\begin{center}
by
\end{center}

\begin{center}
\textsc{Vladimir Lebedev} (Moscow)
\end{center}

\begin{quotation}
{\small \textbf{Abstract.} We consider sets in the real line that
have Littlewood--Paley properties $\mathrm{LP}(p)$ or
$\mathrm{LP}$ and study the following question: How thick can
these sets be?

  References: 14 items.


  }
\end{quotation}

\quad

   Let $E$ be a closed Lebesgue measure zero set
in the real line $\mathbb R$ and let $I_k, ~k=1, 2, \ldots,$ be
the intervals complimentary to $E$, i.e., the connected components
of the compliment $\mathbb R\setminus E$. Call $S_k$ the operator
defined by
$$
\widehat{(S_k f)}=1_{I_k}\cdot \widehat{f},
\qquad f\in L^2\cap L^p(\mathbb R),
$$
where $1_{I_k}$ is the characteristic function of $I_k$, and
$~\widehat{}~$ stands for the Fourier transform. Consider the
corresponding quadratic Littlewood--Paley function:
$$
S(f)=\bigg(\sum_k |S_k f|^2\bigg)^{1/2}.
$$
Following [12] we say that $E$ has property $\mathrm{LP}(p)
~(1<p<\infty)$ if for all $f\in L^p(\mathbb R)$ we have
$$
c_1\|f\|_{L^p(\mathbb R)}\leq\|S(f)\|_{L^p(\mathbb R)}\leq
c_2 \|f\|_{L^p(\mathbb R)},
$$
where $c_1, c_2$ are positive constants independent of $f$. In the
case when a set has property $\mathrm{LP}(p)$ for all $p,
~1<p<\infty,$ we say that it has property $\mathrm{LP}$.

  The role of such sets in harmonic analysis and particularly
in multiplier theory is well-known. We recall that if $G$ is a
locally compact Abelian group and $\Gamma$ is the group dual to
$G$, then a function $m \in L^\infty (\Gamma)$ is called an $L^p$
-Fourier multiplier, $1\leq p\leq \infty$, if the operator $Q$
given by
$$
\widehat{Qf}=m\cdot\widehat{f}, \qquad f\in L^p\cap
L^2(G),
$$
is a bounded operator from $L^p(G)$ to itself (here $~\widehat{}~$
is the Fourier transform on $G$). The space of all these
multipliers is denoted by $M_p(\Gamma)$. Provided with the norm
$$
\|m\|_{M_p(\Gamma)}=\|Q\|_{L^p(G)\rightarrow L^p(G)},
$$
the space $M_p(\Gamma)$ is a Banach algebra (with the usual
multiplication of functions). For basic facts on multipliers in
the cases when $\Gamma=\mathbb R, ~\mathbb Z, ~\mathbb T$, where
$\mathbb Z$ is the group of integers and $\mathbb T=\mathbb
R/2\pi\mathbb Z$ is the circle, see [1], [13, Chap. IV], [7].

  A classical example of an infinite set that has property
$\mathrm{LP}$ is the set $E=\{\pm 2^k, ~k\in \mathbb Z \}\cup\{0\}$ (see,
e.g., [13, Chap. IV, Sec. 5]). From arithmetic and combinatorial point of
view sets that have property $\mathrm{LP}(p)$ or $\mathrm{LP}$ were studied
extensively, see, e. g., the works [1]--[3], [12]. With the exclusion of [12]
these works deal with countable sets, particularly, with subsets of $\mathbb
Z$. At the same time there exist uncountable sets that have property
$\mathrm{LP}$. This fact was first established by Hare and Klemes in [3], see
also [8] and [9, Sec. 4].

   In this paper we study the following question: How
thick can a set $E\subseteq\mathbb R$ that has property $\mathrm{LP}(p)$
($p\neq 2$) or property $\mathrm{LP}$ be? In Theorems 1 and 2 we show that
such a set can not be metrically very thick, namely it is porous and the
measure of the $\delta$ -neighbourhood of any portion of it tends to zero
quite rapidly (as $\delta\rightarrow+0$). As a consequence we obtain (see
Corollary) an estimate for the Hausdorff dimension of these sets. An
immediate consequence of our estimate is that if a set has property
$\mathrm{LP}$, then its Hausdorff dimension is equal to zero. In Theorem 3 we
show that there exist sets which are thin in several senses simultaneously
but have $\mathrm{LP}(p)$ property for no $p\neq 2$. In Theorem 4 we show
that a set can be quite thick but at the same time have property
$\mathrm{LP}$. In part our arguments are close to those used by other authors
to study subsets of $\mathbb Z$ but the mere fact of existence of uncountable
(i.e. thick in the sense of cardinality) sets that have property
$\mathrm{LP}$ brings some specific details to the subject.

   It is well-known that a set has property $\mathrm{LP}(p)$
if and only if it has property $\mathrm{LP}(q)$, where $1/p+1/q=1$
(see, e.g., [12]). Thus, it suffices to consider the case when
$1<p<2$.

   We use the following notation. For a set
$F\subseteq\mathbb R$ we denote its open $\delta$ -neighbourhood
($\delta>0$) by $(F)_\delta$. If $F$ is measurable, then $|F|$
means its Lebesgue measure. A portion of a set $F\subseteq\mathbb
R$ is a set of the form $F\cap I$, where $I$ is a bounded
interval. By $\dim F$ we denote the Hausdorff dimension of $F$.
For basic properties of the Hausdorff dimension we refer the
reader to [11]. For a set $F\subseteq\mathbb R$ and a point
$t\in\mathbb R$ we put $F+t=\{x+t: x\in F\}$. By
$\mathrm{card}\,A$ we denote the number of elements of a finite
set $A$. By arithmetic progression of length $N$ we mean a set of
the form $\{a+kd, ~k=1, 2, \ldots, N\}$, where $a, ~d\in\mathbb R$
and $d\neq 0$. We use $c, ~c(p), ~c(p, E)$... to denote various
positive constants which may depend only on $p$ and the set $E$.

   We recall that a set $F\subseteq\mathbb R$ is said to be porous
if there exists a constant $c>0$ such that every bounded interval
$I\subseteq\mathbb R$ contains a subinterval $J$ with $|J|\geq c
|I|$ and $J\cap F=\varnothing$.

\quad

\textsc{Theorem 1.} \emph{Let $E\subseteq\mathbb R$ be a closed
set of measure zero. Suppose that $E$ has property
$\mathrm{LP}(p)$ for some $p, ~p\neq 2$. Then $E$ is porous.}

\quad

   Earlier Hare and Klemes showed that if a set
in $\mathbb Z$ has property $\mathrm{LP}$ then it is porous [2,
Theorem 3.7].

  To prove Theorem 1 we need certain lemmas.

\quad

\textsc{Lemma 1.} \emph{Let $1<p<\infty$. Let $\varphi : \mathbb
R^n\rightarrow\mathbb R$ be a nonconstant affine mapping. Suppose
that a function $m\in M_p(\mathbb R)$ is continuous at each point
of the set $\varphi(\mathbb Z^n)$. Then the restriction
$m\circ\varphi_{|\mathbb Z^n}$ of the superposition
$m\circ\varphi$ to $\mathbb Z^n$ belongs to $M_p(\mathbb Z^n)$,
and $\|m\circ\varphi_{|\mathbb Z^n}\|_{M_p(\mathbb Z^n)}\leq c
\|m\|_{M_p(\mathbb R)}$, where $c=c(p)>0$ is independent of
$\varphi, ~m$ and the dimension $n$.}

\quad

\emph{Proof.} The proof is a trivial combination of the two
well-known assertions on multipliers. The first one is the theorem
on superpositions with affine mappings [4, Chap. I, Sec. 1.3 ],
which implies that for every function  $m\in M_p(\mathbb R)$ we
have $m\circ\varphi\in M_p(\mathbb R^n)$ and
$\|m\circ\varphi\|_{M_p(\mathbb R^n)}=\|m\|_{M_p(\mathbb R)}$. The
second one is the de Leeuw theorem [10] (see also [5]) on
restrictions to $\mathbb Z^n$, according to which if a function
$g\in M_p(\mathbb R^n)$ is continuous at the points of $\mathbb
Z^n$, then its restriction $g_{|\mathbb Z^n}$ to $\mathbb Z^n$
belongs to $M_p(\mathbb Z^n)$ and $\|g_{|\mathbb
Z^n}\|_{M_p(\mathbb Z^n)}\leq c(p)\|g\|_{M_p(\mathbb R^n)}$. The
lemma is proved.

\quad

\textsc{Lemma 2.} \emph{Let $E\subseteq\mathbb R$ be a nowhere
dense set and let $F\subseteq\mathbb R$ be a finite or countable
set. Then for each $\delta >0$ there exists $\xi\in\mathbb R$ such
that $|\xi|<\delta$ and $(F+\xi)\cap E=\varnothing$.}

\quad

\emph{Proof.} The set
$$
\bigcup_{t \in F}(E-t),
$$
being a union of at most countable family of nowhere dense sets,
can not contain the whole interval $(-\delta, ~\delta)$, hence
there exists $\xi\in(-\delta, ~\delta)$ that does not belong to
the union. The lemma is proved.

\quad

  We say that a (finite or countable) set $F\subseteq\mathbb R$
splits a closed set $E\subseteq\mathbb R$ if $F\subseteq\mathbb R\setminus E$ and there are no two distinct points of $F$ contained in the same interval
complimentary to $E$.

\quad

\textsc{Lemma 3.} \emph{Let $1<p<2$. Let $E\subseteq\mathbb R$ be a set that
has property $\mathrm{LP}(p)$. Suppose that $F$ is a subset of an arithmetic
progression of length $N$ and suppose that $F$ splits $E$. Then
$\mathrm{card}\,F\leq c(p, E) N^{2/q}$, where $1/p+1/q=1$.}

\quad

\emph{Proof.} This lemma can be deduced from Theorems 1.2 and 1.3 of the work [12]. We give an independent short and simple proof based on quite standard argument. Consider an arithmetic progression $\{a+kd, ~k=1, 2, \ldots, N\}$. We can assume that $d>0$. Suppose that a set $F=\{a+k_jd, ~j=1, 2, \ldots, \nu\}$, where $1\leq k_j\leq N$, splits $E$. For $j=1, 2, \ldots, \nu$ let $\Delta_j$ be the interval of length $\delta$ centered at $a+k_jd$, where $\delta>0$ is so small that $\delta<d$ and $\Delta_j\cap E=\varnothing, ~j=1, 2, \ldots, \nu$. We put
$$
m_\theta =\sum_{j=1}^\nu r_j(\theta)\cdot 1_{\Delta_j},
$$
where $r_j(\theta)=\textrm{sign} \sin 2^j\pi \theta, ~\theta\in
[0,1], ~j=1, 2, \ldots,$ are the Rademacher functions.

It is well-known that if a set $E$ has property $\mathrm{LP}(p)$, then it has Marcinkiewicz property $\mathrm{Mar}(p)$, namely \footnote{Actually properties $\mathrm{LP(p)}$ and $\mathrm{Mar(p)}$ are equivalent, see, e.g., [12, Theorem 1.1].}, for each function $m\in L^\infty(\mathbb R)$, whose variations $\textrm{Var}_{I_k}~m$ on the intervals $I_k$ complimentary to $E$ are uniformly bounded, we have $m\in M_p(\mathbb R)$ and
$$
\|m\|_{M_p(\mathbb R)}\leq c(p, E)\Big(\|m\|_{L^\infty (\mathbb
R)}+\sup_k \textrm{Var}_{I_k}~m \Big).
\eqno(1)
$$
Thus we have $\|m_\theta\|_{M_p(\mathbb R)}\leq c$, where $c>0$ is
independent of $N$ and $\theta$. Consider the affine mapping
$\varphi(x)=a+dx, ~x\in\mathbb R$. Using Lemma 1 in the case when
$n=1$, we see that
$$
\|m_{\theta}\circ\varphi_{| \mathbb
Z}\|_{M_p(\mathbb Z)}\leq c(p)\|m_\theta\|_{M_p(\mathbb R)}\leq
c_1(p).
$$
Thus
$$
\bigg\|\sum_k m_\theta(a+kd)c_ke^{ikx}\bigg\|_{L^p(\mathbb T)}\leq
c_1(p)\bigg\|\sum_k c_ke^{ikx}\bigg\|_{L^p(\mathbb T)}
$$
for every trigonometric polynomial $\sum_k c_ke^{ikx}$. In
particular,
$$
\bigg\|\sum_{k=1}^N m_\theta(a+kd)e^{ikx}\bigg\|_{L^p(\mathbb T)}\leq
c_1(p)\bigg\|\sum_{k=1}^N e^{ikx}\bigg\|_{L^p(\mathbb T)}.
$$
Hence,
$$
\bigg\|\sum_{j=1}^\nu r_j(\theta)e^{ik_{j}x}\bigg\|_{L^p(\mathbb
T)}\leq c_1(p)\bigg\|\sum_{k=1}^N e^{ikx}\bigg\|_{L^p(\mathbb T)}.
\eqno(2)
$$

  It is easy to verify that
$$
\bigg\|\sum_{k=1}^N e^{ikx}\bigg\|_{L^p(\mathbb T)}\leq c(p) N^{1/q},
$$
so, (2) yields
$$
\int_\mathbb T \bigg|\sum_{j=1}^\nu
r_j(\theta)e^{ik_{j}x}\bigg|^p~dx\leq c_2(p) N^{p/q}.
$$
By integrating this inequality with respect to $\theta\in [0,1]$
and using the Khintchine inequality:
$$
\Bigg (\int_0^1 \bigg |\sum_j c_jr_j(\theta) \bigg |^p d\theta
\Bigg )^{1/p}\geq c\bigg(\sum_j|c_j|^2\bigg)^{1/2}, \qquad
1\leq p<2,
$$
(see, e.g., [14, Chap. V, Sec. 8]), we obtain $\nu^{p/2}\leq
c_3(p) N^{p/q}$. The lemma is proved.

\quad

\emph{Proof of Theorem 1.} We can assume that $1<p<2$. For a
bounded interval $I\subseteq\mathbb R$ let
$$
d(I)=\sup\{|J|: J\textrm{ is an interval}, ~J\subseteq I, ~J\cap
E=\varnothing\}.
$$
Suppose that $E$ is not porous. Then, for each positive integer
$N$ we can find a (bounded) interval $I$ such that
$0<d(I)<|I|/3N$. Let $d=2d(I)$. Consider an arithmetic progression
$t_k=a+kd, ~k=1, 2, \ldots, N,$ that lies in the interior of $I$.
Using Lemma 2, we can find $\xi$ such that $t_k+\xi \notin E,
~k=1, 2, \ldots, N,$ and $\xi$ is so small that $\{t_k+\xi, ~k=1,
2, \ldots, N\}\subseteq I$. Note that since $d=2d(I)$, there no
two distinct points of the progression $\{t_k+\xi, ~k=1, 2,
\ldots, N\}$ that lie in the same interval complimentary to $E$.
Thus this progression splits $E$. By Lemma 3 this is impossible if
$N$ is sufficiently large. The theorem is proved.

\quad

\textsc{Theorem 2.} \emph{Let $1<p<2$. Let $E\subseteq\mathbb R$
be a closed set of measure zero. Suppose that $E$ has property
$\mathrm{LP}(p)$. Then each portion $E\cap I$ of $E$ satisfies
$$
|(E\cap I)_\delta|\leq c|I|^{2/q}\delta^{1-2/q},
$$
where $1/p+1/q=1$ and the constant $c=c(p, E)>0$ is independent of
$I$ and $\delta$.}

\quad

   Theorem 2 immediately implies an estimate for the
Hausdorff dimension of sets that have $\mathrm{LP}(p)$ property.
Namely, the following corollary is true.

\quad

\textsc{Corollary.} \emph{If $1<p<2$ and a set $E\subseteq\mathbb
R$ has property $\mathrm{LP}(p)$, then $\dim E\leq 2/q$, where
$1/p+1/q=1$. Thus, if  $E$ has property $\mathrm{LP}$, then $\dim
E=0$.}

\quad

\emph{Proof of Theorem 2.} Consider an arbitrary portion $E\cap I$ of the set $E$. Let $J$ be the interval concentric with $I$ and of two times larger length. Denote the left-hand endpoint of $J$ by $a$. Fix a positive integer $N$ and consider the progression $a+kd, ~k=1, 2, \ldots, N,$ where $d=|J|/N$. By Lemma 2 one can find $\xi$ such that none of the elements of the progression $\{a+kd+\xi,~k=1, 2, \ldots,
N\}$ is contained in $E$ and $I\subseteq J+\xi=(a+\xi, ~a+Nd+\xi)$.

   We define intervals $J_k$ by
$$
J_k=(a+(k-1)d+\xi, ~a+kd+\xi), \qquad k=1, 2, \ldots, N.
$$
Consider the intervals $J_{k_j}$ such that $J_{k_j}\cap
E\neq\varnothing$. Obviously their right-hand endpoints split $E$,
so, by Lemma 3, their number is at most $c(p) N^{2/q}$. Thus the
set $E\cap I$ is covered by at most $c(p) N^{2/q}$ intervals of
length $d=2|I|/N$.

  Let $\delta>0$. We can assume that $\delta<|I|$
(otherwise the assertion of the theorem is trivial). Choosing a
positive integer $N$ so that
$$
\frac{2|I|}{N}\leq\frac{\delta}{3}<\frac{4|I|}{N},
$$
we see that the portion $E\cap I$ can be covered by at most
$c(p)(12|I|/\delta)^{2/q}$ intervals of length $\delta/3$. It
remains to replace each of these intervals with the corresponding
concentric interval of nine times larger length. The theorem is
proved. The corollary follows.

\quad

  We note now that a set can be quite thin and at the same
time have property $\mathrm{LP}(p)$ for no $p\neq 2$. Consider a
set
$$
F=\bigg\{\sum_{k=1}^\infty \varepsilon_k l_k, ~\varepsilon_k=0
~\textrm{or}~ 1\bigg \},
\eqno(3)
$$
where $l_k, ~k=1, 2, \ldots,$ are positive numbers with
$l_{k+1}<l_k/2$. It was shown by Sj\"{o}gren and Sj\"{o}lin [12]
that such sets have property $\mathrm{LP}(p)$ for no $p, ~p\neq
2$. (In particular, the Cantor triadic set does not have property
$\mathrm{LP}(p)$ for $p\neq 2$.) Taking a rapidly decreasing
sequence $\{l_k\}$ one can obtain a set $F$ of the form (3) such
that it is porous and the measure of its $\delta$ -neighbourhood
rapidly tends to zero. Still, in a sense, any set of the form (3)
is thick, it is uncountable and all its points are its
accumulation points. Theorem 3 below shows that a set can be thin
in several senses simultaneously, and at the same time have
property $\mathrm{LP}(p)$ for no $p, ~p\neq 2$.

\quad

\textsc{Theorem 3.} \emph{Let $\psi$ be a positive function on an
interval $(0, \delta_0), ~\delta_0>0,$ with
$\lim_{\delta\rightarrow+0}\psi(\delta)/\delta=+\infty$. There
exists a strictly increasing bounded sequence $a_1<a_2<\ldots$
such that the set
$E=\{a_k\}_{k=1}^\infty\cup\{\lim_{k\rightarrow\infty}a_k\}$
satisfies the following conditions: \emph{1)} $E$ is porous;
\emph{2)} $|(E)_\delta|\leq\psi(\delta)$ for all sufficiently
small $\delta>0$; \emph{3)} $E$ has property $\mathrm{LP}(p)$ for
no $p, ~p\neq 2$.}

\quad

\emph{Proof.} Given (real) numbers $a, ~l_1, l_2,\ldots, l_n$
consider the set of all points $a+\sum_{j=1}^n \varepsilon_j~l_j$,
where $\varepsilon_j=0 ~\textrm{or}~ 1$. Assume that the
cardinality of this set is $2^n$. Following [6] we call such a set
an $n$ -chain.\footnote{An $n$ -chain is a particular case of what is called a parallelepiped of dimension $n$, that is of a sum of $n$ two-element sets.}

  We shall need the following refinement of the
Sj\"{o}gren and Sj\"{o}lin result on the sets (3). This refinement also provides a partial extention of Proposition 3.4 of the work [2], that treats subsets of integers, to the general case of closed measure zero sets in the line.

\quad

  \textsc{Lemma 4.} \emph{Let $E\subseteq\mathbb R$ be a closed
set of measure zero. Suppose that for an arbitrary large $n$ the
set $E$ contains an $n$ -chain. Then $E$ has property
$\mathrm{LP}(p)$ for no $p\neq 2$.}

\quad

\emph{Proof.} Suppose that, contrary to the assertion of the
lemma, $E$ has property $\mathrm{LP}(p)$ for some $p, ~p\neq 2$.
We can assume that $1<p<2$.

   Let $n$ be such that $E$ contains an $n$ -chain
$$
a+\sum_{j=1}^n \varepsilon_j~l_j,  \quad
(\varepsilon_1, \varepsilon_2, \ldots, \varepsilon_n)\in\{0;1\}^n.
\eqno(4)
$$

   Consider the set
$$
B=\bigg\{a+\sum_{j=1}^n k_j~l_j,
~(k_1, k_2, \ldots, k_n)\in\mathbb Z^n\bigg\}.
$$
By Lemma 2 there exists an arbitrary small $\xi$ such that
$$
(B+\xi)\cap E=\varnothing
\eqno(5)
$$
Clearly, if $\xi$ is small enough, then no two distinct points of
the chain obtained by the same shift $\xi$ of the chain (4) can
lie in the same interval complimentary to $E$. Thus, there exists
$\xi$ such that (5) holds and simultaneously the $n$ -chain
$$
a+\xi+\sum_{j=1}^n \varepsilon_j~l_j,  \quad
(\varepsilon_1, \varepsilon_2, \ldots, \varepsilon_n)\in\{0;1\}^n,
$$
splits $E$.

   For each $\varepsilon=
(\varepsilon_1, \varepsilon_2, \ldots, \varepsilon_n)\in\{0;1\}^n$
let $I_\varepsilon$ denote the interval complimentary to $E$ that
contains the point $a+\xi+\sum_{j=1}^n \varepsilon_j~l_j$. For an
arbitrary choice of signs $\pm$ consider the function
$$
m=\sum_{\varepsilon\in\{0;1\}^n }\pm 1_{I_\varepsilon}.
$$
We have (see (1))
$$
\|m\|_{M_p(\mathbb R)}\leq c,
\eqno(6)
$$
where $c>0$ is independent of $n$ and the choice of signs.

  Consider the following affine mapping $\varphi$:
$$
\varphi(x)=a+\xi+
\sum_{j=1}^n x_j~l_j,
\quad x=(x_1, x_2, \ldots, x_n)\in\mathbb R^n.
$$

  Note that condition (5) implies that the function $m$ is continuous
at each point of the set $\varphi(\mathbb Z^n)$. Using Lemma 1, we
obtain (see (6)) $m\circ\varphi_{|\mathbb Z^n}\in M_p(\mathbb
Z^n)$ and
$$
\|m\circ \varphi_{|\mathbb Z^n}\|_{M_p(\mathbb Z^n)}\leq c,
$$
where the constant $c>0$ is independent of $n$ and the choice of
signs.

  Therefore, for an arbitrary trigonometric polynomial
$$
\sum_{k\in \mathbb Z^n}c_ke^{i(k,t)}
$$
on the torus $\mathbb T^n$ we have
$$
\bigg\|\sum_{k\in \mathbb Z^n}m\circ \varphi
(k)c_ke^{i(k,t)}\bigg\|_{L^p(\mathbb T^n)}\leq
c\bigg\|\sum_{k\in\mathbb Z^n}c_ke^{i(k,t)}\bigg\|_{L^p(\mathbb
T^n)}.
$$
(We use $(k, t)$ to denote the usual inner product of vectors
$k\in\mathbb Z^n$ and $t\in\mathbb T^n$.) In particular, taking
$c_k=1$ for $k\in\{0;1\}^n$ and $c_k=0$ for $k\notin\{0;1\}^n$, we
obtain
$$
\bigg\|\sum_{\varepsilon\in\{0;1\}^n
}m\Big(a+\xi+\sum_{j=1}^n\varepsilon_j
l_j\Big)e^{i(\varepsilon,t)}\bigg\|_{L^p(\mathbb T^n )}\leq
c \bigg\|\sum_{\varepsilon\in\{0;1\}^n}e^{i(\varepsilon,t)}
\bigg\|_{L^p(\mathbb
T^n )}.
$$
That is
$$
\bigg\|\sum_{\varepsilon\in\{0;1\}^n }\pm
e^{i(\varepsilon,t)}\bigg\|_{L^p(\mathbb T^n )}\leq c
\bigg\|\sum_{\varepsilon\in\{0;1\}^n}e^{i(\varepsilon,t)}
\bigg\|_{L^p(\mathbb
T^n )}.
$$

   Raising this inequality to the power $p$ and averaging
with respect to the signs $\pm$ (i.e., using the Khintchine
inequality), we obtain
$$
\bigg\|\sum_{\varepsilon\in\{0;1\}^n}e^{i(\varepsilon,t)}
\bigg\|_{L^2(\mathbb
T^n )}\leq c
\bigg\|\sum_{\varepsilon\in\{0;1\}^n}e^{i(\varepsilon,t)}
\bigg\|_{L^p(\mathbb
T^n )}. \eqno(7)
$$

  Note that
$$
\sum_{\varepsilon\in\{0;1\}^n}e^{i(\varepsilon,t)}=\prod_{j=1}^n
(1+e^{it_j}), \quad t=(t_1, t_2, \ldots, t_n)\in\mathbb T^n,
$$
so (7) yields
$$
\|1+e^{it}\|_{L^2(\mathbb T)}^n\leq
c\|1+e^{it}\|_{L^p(\mathbb T)}^n.
\eqno(8)
$$
Since $n$ can be arbitrarily large, relation (8) implies
$$
\|1+e^{it}\|_{L^2(\mathbb T)}\leq \|1+e^{it}\|_{L^p(\mathbb T)},
$$
which, as one can easily verify, is impossible for $1<p<2$. The
lemma is proved.

\quad

  \textsc{Lemma 5.} \emph{Let $l_k, ~k=1, 2, \ldots,$ be positive
numbers satisfying $l_{k+1}<l_k/2$. Then the set $F$ defined by
\emph{(3)} contains a strictly increasing sequence
$S=\{a_k\}_{k=1}^\infty$ such that for every $n$ the sequence $S$
contains an $n$ -chain.}

\quad

\emph{Proof.} For $n=1, 2, \ldots$ let
$$
\alpha_n=\sum_{k=1}^{n^2} l_k, \qquad\beta_n=\sum_{k=1}^{n^2+n} l_k.
$$
Clearly $\alpha_1<\beta_1<\alpha_2<\beta_2<\ldots$, so the closed
intervals $[\alpha_n, \beta_n], ~n=1, 2, \ldots,$ are pairwise
disjoint.

  Define sets $F_n\subseteq F, ~n=1, 2, \ldots,$ as follows
$$
F_n=\bigg\{l_1+l_2+\ldots+l_{n^2}
+\sum_{k=n^2+1}^{n^2+n}\varepsilon_k l_k,
~\varepsilon_k=0 ~\textrm{or}~1\bigg\}.
$$
Note that $F_n\subseteq [\alpha_n, \beta_n]$ for all $n=1,
2,\ldots$.

  It remains to put
$$
S=\bigcup_{n=1}^\infty F_n.
$$
The lemma is proved.

\quad

    We shall now complete the proof of the theorem.
Replacing, if needed, the function $\psi(\delta)$ with
$$
\widetilde{\psi}(\delta)=\delta\inf_{0<t\leq\delta}\frac{\psi(t)}{t},
$$
we can assume that the relation $\psi(\delta)/\delta$ increases to
$+\infty$ as $\delta$ decreases to zero.

   Take a strictly increasing sequence of positive integers
$n_k, ~k=1, 2, \ldots,$ so that
$$
6\cdot 2^k\leq\frac{\psi(3^{-n_k})}{3^{-n_k}}, \qquad k=1, 2, \ldots.
\eqno(9)
$$
Consider the set
$$
F=\Big\{\sum_{k=1}^\infty \varepsilon_k 3^{-n_k}, ~\varepsilon_k=0
~\textrm{or}~ 1\Big \}.
$$
It is clear that $F$ is porous (as a subset of the Cantor triadic
set).

  Assuming that $\delta>0$ is sufficiently small, we can find
$k$ such that
$$
3^{-n_{k+1}}\leq \delta <3^{-n_k}.
\eqno(10)
$$
Note that $F$ can be covered by $2^{k+1}$ closed intervals of
length $3^{-n_{k+1}}$. Consider the $\delta$ -neighbourhood of
each of these intervals. We see that (see (10))
$$
|(F)_\delta|\leq 2^{k+1} 3 \delta.
$$
Hence, taking (9), (10) into account, we obtain
$$
|(F)_\delta|\leq\frac{\psi(3^{-n_k})}{3^{-n_k}}\delta\leq \psi(\delta).
$$

  Using Lemma 5 we can find a strictly increasing sequence
$S=\{a_k\}_{k=1}^\infty$ contained in $F$, such that for every $n$
the sequence $S$ contains an $n$ -chain. Let $E=S\cup\{a\}$, where
$a=\lim_{k\rightarrow\infty}a_k$. It remains to use Lemma 4. The
theorem is proved.

\quad

  Our next goal is to construct a set that has property
$\mathrm{LP}(p)$ or property $\mathrm{LP}$ and at the same time is thick.
Theorem 2 implies that if $1<p<2$ and a bounded set $E$ has property
$\mathrm{LP}(p)$, then $|(E)_\delta|=O(\delta^{1-2/q})$ as
$\delta\rightarrow+0$. Hence, if a bounded set $E$ has property
$\mathrm{LP}$, then $|(E)_\delta|=O(\delta^{1-\varepsilon})$ for all
$\varepsilon>0$. The author does not know if these estimates are sharp. A
partial solution to this problem is given by Theorem 4 below. This theorem is
a simple consequence of the Hare and Klemes theorem [3, Theorem A], which
provides a sufficient condition for a set to have property $\mathrm{LP}(p)$.
Stated for sets in $\mathbb Z$ this theorem, as is noted at the end of the
work [3], easily transfers to sets in $\mathbb R$ and allows to construct
perfect sets that have this property. We shall use the version of the Hare
and Klemes theorem stated in [9, Sec. 4]. According to this version, for each
$p, ~1<p<\infty,$ there is a constant $\tau_p ~(0<\tau_p<1)$ with the
following property. Let $E$ be a closed set of measure zero in the interval
$[0, 1]$. Suppose that, under an appropriate enumeration, the intervals $I_k,
~k=1, 2, \ldots,$ complimentary to $E$ in $[0, 1]$ (i.e., the connected
components of the compliment $[0, 1]\setminus E$) satisfy
$$
\delta_{k+1}/\delta_k\leq \tau_p, \quad k=1, 2, \ldots,
\eqno(11)
$$
where $\delta_k=|I_k|$. Then $E$ has property $\mathrm{LP}(p)$.
This in turn implies that if
$$
\lim_{k\rightarrow\infty}\delta_{k+1}/\delta_k=0,
\eqno(12)
$$
then $E$ has property $\mathrm{LP}$.

\quad

\textsc{Theorem 4.}

(a) \emph{Let $1<p<\infty$. There exists a perfect set $E\subseteq
[0, 1]$ which has property $\mathrm{LP}(p)$ and at the same
time $|(E)_\delta|\geq c\,\delta\log 1/\delta$ for all
sufficiently small $\delta>0$.}

(b) \emph{Let $\gamma(\delta)$ be a positive non-decreasing
function on $(0, +\infty)$ with
$\lim_{\delta\rightarrow+0}\gamma(\delta)=0$. There exists a
perfect set $E\subseteq [0, 1]$ which has property
$\mathrm{LP}$ and at the same time $|(E)_\delta|\geq
c\,\gamma(\delta)\delta\log 1/\delta$.}

\quad

\emph{Proof.} Let $\delta_k, ~k=1, 2, \ldots,$ be a sequence of
positive numbers with
$$
\sum_k \delta_k=1.
\eqno(13)
$$
Let $E\subseteq [0, 1]$ be a closed set. Assume that, under an
appropriate enumeration, the intervals $I_k, ~k=1, 2, \ldots,$
complimentary to $E$ in $[0, 1]$ satisfy $|I_k|=\delta_k, ~k=1, 2,
\ldots$. In this case we say that $E$ is generated by the sequence
$\{\delta_k\}$. (Certainly $|E|=0$.) Note that for each sequence
$\{\delta_k\}$ of positive numbers with (13) there exists a
perfect set $E\subseteq [0, 1]$ generated by $\{\delta_k\}$.

  It is easy to see that if $E$ is a set generated by a positive
sequence $\{\delta_k\}$ satisfying (13), then for all $\delta>0$
we have
$$
|(E)_\delta|\geq  2\delta\, \mathrm{card}\{k : \delta_k>2\delta\}.
\eqno(14)
$$
Indeed, if $I_k=(a_k, b_k)$ is an arbitrary interval complimentary
to $E$ in $[0, 1]$ such that $|I_k|>2\delta$, then the $\delta$
-neighbourhood of $E$ contains the intervals $(a_k, a_k+\delta)$
and $(b_k-\delta, b_k)$.

  We shall prove part (a) of the theorem. Fix $p, ~1<p<\infty$.
Let
$$
\delta_k=a e^{-kb}, \quad k=1, 2, \ldots,
$$
where the positive constants $a$ and $b$ are chosen so that
conditions (11), (13) hold. Consider a perfect set $E\subseteq
[0,1]$ generated by the sequence $\{\delta_k\}$. Using estimate
(14), we see that
$$
|(E)_\delta|\geq 2\delta \bigg (\frac{1}{b} \log\frac{a}{2\delta}-1\bigg),
$$
which proves (a).

  Now we shall prove (b). Without loss of generality
we can assume that $\gamma(1/e)=1/4$. Let
$$
b(x)=\frac{1}{\gamma(e^{-x})}, \quad x>0.
$$
The function $b$ is non-decreasing, $b(x)\rightarrow+\infty$ as
$x\rightarrow +\infty$, and $b(1)=4$.

    Define a sequence $\{\delta_k\}$ by
$$
\delta_k=a e^{-kb(k)}, \quad k=1, 2, \ldots,
$$
where the constant $a>0$ is chosen so that condition (13) holds.

   Note that
$$
\delta_{k+1}/\delta_k=e^{-((k+1)b(k+1)-kb(k))}
\leq e^{-b(k)}\rightarrow 0, \qquad k\rightarrow\infty,
$$
and thus, (12) holds.

  Consider a perfect set $E\subseteq [0, 1]$ generated by the sequence
$\{\delta_k\}$.

  Let $\delta>0$ be sufficiently small. Chose a positive integer
$k=k(\delta)$ so that
$$
\delta_{k+1}\leq 2\delta<\delta_k.
\eqno(15)
$$
We have
$$
\mathrm{card}\{k : \delta_k>2\delta\}\geq k(\delta).
$$
So (see (14)),
$$
|(E_\delta)|\geq 2\delta k(\delta).
\eqno(16)
$$

   Note that (15) implies
$$
kb(k)<\log\frac{a}{2\delta}\leq (k+1)b(k+1).
$$
Hence, for all sufficiently small $\delta>0$ we have
$$
\frac{1}{2}kb(k)<\log\frac{1}{\delta}\leq 2(k+1)b(k+1).
\eqno(17)
$$
The left-hand inequality in (17) yields (recall that $b(1)=4$)
$$
2k=\frac{1}{2}k b(1)\leq\frac{1}{2}kb(k)<\log\frac{1}{\delta},
$$
whence
$$
b(2k)\leq b\bigg(\log\frac{1}{\delta}\bigg)=\frac{1}{\gamma(\delta)}.
$$
Combining this inequality and the right-hand inequality in (17),
we see that
$$
\log\frac{1}{\delta}\leq 2(k+1)b(k+1)\leq 4kb(2k)
\leq 4k\frac{1}{\gamma(\delta)}.
$$
So,
$$
\frac{1}{4}\gamma(\delta)\log\frac{1}{\delta}\leq k=k(\delta).
$$

   Thus (see (16)),
$$
|(E)_\delta|\geq \frac{1}{2}\gamma(\delta)\delta\log\frac{1}{\delta}.
$$
The theorem is proved.

\quad

\quad

\textsc{Remark.}  As far as the author knows, the question on the
existence of a set that has property $\mathrm{LP}(p)$ for some $p,
~p\neq 2,$ but does not have property $\mathrm{LP}$ is open.

\begin{center}
\textbf{References}
\end{center}

\flushleft
\begin{enumerate}

\item R. E. Edwards and  G. I. Gaudry, \emph{Littlewood-Paley
    and multiplier theory}, Springer–Verlag,
    Berlin–Heidelberg–New York, 1977.

\item K. E. Hare and I. Klemes, \emph{Properties of
    Littlewood-Paley sets}, Math. Proc. Camb. Phil. Soc.,
    \textbf{105} (1989), 485--494.

\item K. E. Hare and  I. Klemes, \emph{On permutations of
    lacunary intervals}, Trans. Amer. Math. Soc., 347 (1995),
    no. 10, 4105--4127.

\item L. H\"ormander, \emph{Estimates for translation
    invariant operators in $L^p$ spaces}, Acta Math., 104
    (1960), 93–140.

\item M. Jodeit, \emph{Restrictions and extensions of Fourier
    multipliers}, Studia Math., \textbf{34} (1970), 215--226.

\item J.-P. Kahane, \emph{S\'eries de Fourier absolument
    convergentes}, Springer–Verlag, Berlin–Heidelberg–New
    York, 1970.

\item R. Larsen, \emph{An introduction to the theory of
    multipliers}, Springer–Verlag, Berlin–Heidelberg–New York,
    1971.

\item V. Lebedev and  A. Olevski\v{\i}, \emph{Bounded groups
    of translation invariant operators}, C. R. Acad. Sci.
    Paris., S\'{e}rie I, 322 (1996), 143--147.

\item V. V. Lebedev,  A. M. Olevski\v{\i}, \emph{$L^p$
    -Fourier multipliers with bounded powers}, Izv. Math., 70
    (2006), no. 3, 549-585.

\item K. de Leeuw, \emph{On $L^p$ multipliers}, Ann. Math., 81
    (1965), no. 2, 364--379.

\item P. Mattila, \emph{Geometry of sets and measures in
    Euclidean spaces}, Cambrige Univ. Press, 1995.

\item P. Sj\"{o}gren  and  P. Sj\"{o}lin,
    \emph{Littlewood-Paley decompositions and Fourier
    multipliers with singularities on certain sets}, Ann.
    Inst. Fourier, Grenoble, 31 (1981), no. 1, 157--175.

\item E. M. Stein, \emph{Singular integrals and
    differentiability properties of functions}, Princeton
    Math. Ser., vol. 30, Princeton Univ. Press, Princeton, NJ
    1970.

\item A. Zygmund, \emph{Trigonometric series}, vols. I, II,
    Cambridge Univ. Press, New York 1959.

\end{enumerate}

\quad

\quad

National Research University Higher School of Economics,

Moscow Institute of Electronics and Mathematics,

Bolshoi Trekhsvjatitelskii, 3, Moscow, 109028, Russia

e-mail: \emph{lebedevhome@gmail.com}

\end{document}